\documentclass{birkjour}
\usepackage{amssymb}
\usepackage{url}

\usepackage[utf8]{inputenc}
\usepackage[T1]{fontenc}
\usepackage{wasysym}
\usepackage{amsmath}
\usepackage{tikz-cd}
\usepackage{mathrsfs,multicol}
\usepackage[pagebackref, hypertexnames=false, colorlinks, citecolor=red, linkcolor=red]{hyperref}

 \newtheorem{thm}{Theorem}[section]
 
 \newtheorem{lem}[thm]{Lemma}
 \newtheorem{prop}[thm]{Proposition}
 \theoremstyle{definition}
 \newtheorem{defn}[thm]{Definition}
 \theoremstyle{remark}
 \newtheorem{rem}[thm]{Remark}
 
 \numberwithin{equation}{section}

\def\TT{\mathbb{T}}
\def\RR{\mathbb{R}}
\def\ZZ{\mathbb{Z}}
\def\CC{\mathbb{C}}

\def\DD{\mathbb{D}}
\newcommand{\ID}{{\mathbf{1}}}
\newcommand{\al}{{\alpha}}

\newcommand{\f}{{\varphi}}
\newcommand{\IM}{{\operatorname{Im}}}
\newcommand{\RE}{{\operatorname{Re}}}

\newcommand{\bte}{\boldsymbol{\theta}}

\newcommand{\R}{{\mathbb  R}}
\newcommand{\fD}{{\mathcal  D}}

\newcommand{\D}{{\mathbb  D}}
\newcommand{\T}{\mathbb{T}}

\newcommand{\te}{{\theta}}
\newcommand{\Z}{{\mathbb  Z}}

\newcommand{\N}{{\mathbb  N}}
\newcommand{\C}{{\mathbb  C}}

\newcommand{\cK}{{\mathcal K}}

\newcommand{\OZ}{{\mathbf{0}}}

\newcommand{\bA}{{\mathbf A}}
\newcommand{\bU}{{\mathbf U}}
\newcommand{\bB}{{\mathbf B}}

\newcommand{\bJ}{{\mathbf J}}

\newcommand{\balpha}{\boldsymbol{\alpha} }
\newcommand{\bmu}{\boldsymbol{\mu}}
\newcommand{\bnu}{\boldsymbol{\nu}}

\newcommand{\OID}{ I }

\newcommand{\fdot}{\,\cdot\,}

\makeatletter
\def\Ddots{\mathinner{\mkern1mu\raise\p@
\vbox{\kern7\p@\hbox{.}}\mkern2mu
\raise4\p@\hbox{.}\mkern2mu\raise7\p@\hbox{.}\mkern1mu}}
\makeatother

\newcommand{\cH}{\mathcal{H}}
\newcommand{\cC}{\mathcal{C}}

\newcommand{\cA}{\mathcal{A}}
\newcommand{\cE}{ E }

\newcommand{\cD}{\mathcal{D}}

\newcommand{\ciG}{\ci\Gamma}

\DeclareMathOperator{\tr}{tr}

\DeclareMathOperator{\Ker}{Ker}
\DeclareMathOperator{\Dom}{Dom}
\DeclareMathOperator{\rk}{rank}
\DeclareMathOperator{\Ran}{Ran}
\DeclareMathOperator{\spa}{span}
\DeclareMathOperator{\im}{Im}
\DeclareMathOperator{\clos}{clos}

\DeclareMathOperator{\supp}{supp}

\newcommand{\ci}[1]{_{ {}_{\scriptstyle #1}}}
\newcommand{\ti}[1]{_{\scriptstyle \text{\rm #1}}}

\begin{document}

%
%
%
%
%
%
%
%
%

\title[Finite-Rank Perturbations]{Spectral Analysis, Model Theory and Applications of Finite-Rank Perturbations}

\author[Frymark]{Dale~Frymark}
\address{Department of Mathematics\\
Stockholm University\\
Kr\"aftriket 6\\
106 91 Stockholm\\
Sweden}
\email{dale@math.su.se}

\author[Liaw]{Constanze~Liaw}
\address{Department of Mathematical Sciences\\
University of Delware\\
501 Ewing Hall\\
Newark, DE  19716\\
USA\\
and \\
CASPER\\
Baylor University\\
One Bear Place \#97328\\     
 Waco, TX  76798\\
 USA}
\email{liaw@udel.edu}

\thanks{The work of Constanze Liaw was supported by the US National Science Foundation under the grant DMS-1802682.}



\date{\today}
\dedicatory{Dedicated to the memory of R.G.~Douglas,\\a magnificent person, administrator and mathematician.}

\begin{abstract}
This survey focuses on two main types of finite-rank perturbations: self-adjoint and unitary. We describe both classical and more recent spectral results. We pay special attention to singular self-adjoint perturbations and model representations of unitary perturbations.
\end{abstract}

\maketitle
\section{Introduction}\label{s-Intro}

Let $\bA$ be a self-adjoint (possibly unbounded) operator on a separable Hilbert space $\cH$. Fix a $d$-dimensional subspace $\cK\le\cH$. Consider all self-adjoint perturbations $\bA+K$ with $\Ran{K}\subset \cK$. All self-adjoint perturbations $\bA+K$ are formally given by the {\em family of self-adjoint finite-rank perturbations}:
\begin{align}\label{d-bAGamma}
\bA\ciG
=
\bA
+\bB\Gamma\bB^*
\end{align}
for some Hermitian $d\times d$ matrix $\Gamma$, where $\bB:\C^d\to \cK$ is an invertible coordinate operator that takes the standard basis $\{{\bf e}_k\}_{k=1}^d$ of $\C^d$ into a basis $\bB{\bf e}_k$ of $\cK$. 
Reducing our attention to the essence of the problem, we always assume without loss of generality that $\cK$ is cyclic for $\bA$ on $\cH$, that is, $\cH = \clos\spa\{(\bA-z I)^{-1}\cK: z\in \C\setminus\R\}.$
See Section \ref{s-SAFRP} for a more general definition of $\bA\ciG$ which applies when the functions $\bB{\bf e}_k$ do not belong to the Hilbert space $\cH$, but are instead taken from a larger space.

The {\em family of self-adjoint rank-one perturbations} represents a special case of the family of finite-rank perturbations given in equation \eqref{d-bAGamma}, and can be formally given by 
\begin{align}\label{e-rankone}
A_\gamma = A+\gamma (\fdot, \f)\f,
\qquad
\f\in \cK
\end{align}
with parameter $\gamma\in\R$. See Subsection \ref{ss-scales} for the precise definition, as well as Subsection \ref{ss-notation} regarding notation on $\bA$ versus $A$.

Interest in this type of perturbation problem originally arose from the theory of self-adjoint extensions \cite{W}. Natural applications to the variation of boundary conditions of differential operators, in particular Sturm--Liouville operators, were investigated by Aronszajn and Donoghue in the 1950's. Other famous perturbation theoretic results, such as those by von Neumann and Kato--Rosenblum, apply because rank-one perturbations are trace class. The great achievements in this field furnish a rather concrete description of the spectral properties of the perturbed operators $A_\gamma$. See Section \ref{s-rankone} for a sampling of these results.

The spectral theory for quantum mechanical systems (see e.g.~\cite{AK}), large random matrices (see e.g.~\cite{BN}) and free probability probability (see e.g.~\cite{BBCF}), and the decoupling of CMV matrices (see e.g.~\cite[Section 4.5]{SIMOPUC1}) present other standard applications. 
Additional applications to quantum graph theory arise from transforming the graph to a tree by adding \emph{partition vertices} to existing edges and imposing boundary conditions on the partition vertices \cite[Ch.~3]{BerkoKuch}. The number of partition vertices that needs to be added in order to transform a graph into a tree is equal to the first Betti or cyclomatic  number of the original graph, which equals the number of edges minus the number of vertices plus the number of connected components.

In the late 1980's and early 1990's, a surge of interest took place in perturbation theory following the discovery of the celebrated Simon--Wolff criterion, which was used in a proof of Anderson localization for the discrete random Schr\"odinger operator in dimension one. A brief discussion of the Simon--Wolff criterion is included in Subsection \ref{ss-SIMWOL}. 

Given two arbitrary operators on the same Hilbert space, it is generally not easy to find out whether they are related via a rank-one or a finite-rank perturbation. 
The situation is different if we consider two (so-called) Anderson-type Hamiltonians. We refer the reader to Subsection \ref{ss-SIMWOL} for a definition. For now it suffices to know that they are perturbation problems with a random perturbation that is almost surely non-compact. 
Under mild assumptions, the essential part of two realizations of an Anderson-type Hamiltonian are related by a rank-one perturbation (almost surely with respect to the product of the probability measures), see \cite{LiawBJMA}.

\vspace{1em}

Unitary perturbation theory is the other main topic of this survey. Let $\bU$ be a unitary operator on a Hilbert space $\cH$. Fix a $d$-dimensional subspace $\mathcal{R}\le\cH$. Then the set of operators $K$ with $\Ran{K}\subset \mathcal{R}$ that make $\bU+K$ a unitary operator can be parametrized by unitary $d\times d$ matrices. Specifically, there is a bijective coordinate operator $\bJ: \C^d\to \mathcal{R}$ so that $K = \bJ(\alpha - \OID)\bJ^*\bU$ for a unitary $d\times d$ matrix $\alpha$. The created \emph{family of unitary finite-rank perturbations} of $\bU$ is given by
\begin{align}\label{d-ual}
\bU_\alpha = \bU +\bJ(\alpha - \OID)\bJ^*\bU,
\end{align}
with $\alpha$ taken from the unitary $d\times d$ matrices. Without loss of generality, we focus on the domain altered by assuming that $\mathcal{R}$ is a $*$-cyclic subspace for $\bU$, i.e.~we assume that $\cH = \clos\spa \{\bU^{k}\mathcal{R}:k\in \Z\}.$

The special case when $d=1$ is closely related to Aleksandrov--Clark theory, and is described in Subsection \ref{ss-SzNF}. In this setting, the family of perturbations in equation \eqref{d-ual} reduce to the well-known \emph{family of unitary rank-one perturbations}
\begin{align}\label{d-sUalpha}
    U_\al = U+\al (\fdot, U^*\f)\ci\cH \f,
\end{align}
with $\al\in \T$ and $\f\in \mathcal{R}$. Again, see Subsection \ref{ss-notation} for notation.

While self-adjoint and unitary operators are intimately connected via the Cayley transform, it is well-known (see e.g.~\cite[Theorem 4.3.1]{BirmanSol-book_1987}) that this correspondence is not a bijection between the two operator classes. In fact, even when the mappings are well-defined, the Cayley transform does not explicitly take \eqref{d-bAGamma} to its analog \eqref{d-ual}. This can be seen for the rank-one setting in Liaw--Treil \cite[pp.~124--128]{LTSurvey}. Also notice that we encounter some inconveniences arising from unbounded operators in the self-adjoint setting. Of course, the unbounded case is exactly what occurs when dealing with boundary conditions of differential operators and several other applications. The unitary setting, on the other hand, is always restricted to bounded operators (see Remark \ref{r-H-n}).

It is therefore surprising that, in spite of these differences, many results on self-adjoint finite-rank perturbations have analogs in the unitary setting. It is also common to find that the problems raise similar questions, e.g.~about the boundary behavior of analytic functions.

Families of rank-one and finite-rank perturbations seem rather elementary, yet their study has revealed a quite subtle nature. Their complexity is verified by connections to several deep fields of analysis: Nehari interpolation problem, holomorphic composition operators, rigid functions, existence of the limit of the Julia--Carath\'eodory quotient, Carleson embedding, and functional models. Some of these connections are the topic of existing books and surveys, including \cite{cimaross, LTSurvey, poltsara2006, Saksman}.

While writing this survey, it became evident that a complete account of the subject of finite-rank perturbations is worthy of a whole book due to the connections to many other fields of mathematics. We decided to focus on a few aspects, while only briefly mentioning others. For example, some deserving topics such as related function theoretic nuances are not surveyed in detail. We also often refer to existing surveys and books on the topic such as, e.g.~\cite{AK, cimaross, LTSurvey, poltsara2006, Saksman, SIMREV}, in order to not overlap excessively.

It should be noted that some central objects of perturbation theory, such as Aleksandrov Spectral Averaging and Poltoratski's Theorem, appear in the Appendix (Section \ref{s-APP}) for convenience.
\vspace{1em}

Section \ref{s-PRELIM} contains highlights of classical perturbation theory that provide additional context for the more specific results to come. In particular, we focus on aspects of the spectrum that are invariant under different types of perturbations. 

Sections \ref{s-rankone} and \ref{s-introunit} present well-known features of rank-one perturbation theory in the self-adjoint and unitary settings respectively. Section \ref{s-rankone} includes a discussion of singular perturbations, some spectral results (including Aronszajn--Donoghue theory) and Nevanlinna--Herglotz functions, which form the backbone of the theory. 
The unitary setting of Section \ref{s-introunit} is built upon Aleksandrov--Clark theory and features the Sz.-Nagy--Foia\c s and de Branges--Rovnyak approach, as well as the overarching Nikolski--Vasyunin transcription free model theory. The latter reduces to the ones by Sz.-Nagy--Foia\c s and de Branges--Rovnyak by choosing a specific weight. These model representations form rather concrete applications of model theory.

Sections \ref{s-SAFRP} through \ref{s-STUnit} focus on finite-rank perturbations. Where possible, the presentation runs in analogy to Sections \ref{s-rankone} and \ref{s-introunit}.

For finite-rank self-adjoint perturbations the setup (Section \ref{s-SAFRP}) is a bit more involved, and we include information on extension theory, as well as a summary of some mathematical physics applications. In Section \ref{s-SAST} we present known results regarding the spectral analysis of finite-rank perturbations and compare them to Aronszajn--Donoghue theory.

Section \ref{s-UMT} contains information about model spaces culminating in the Nikolski--Vasyunin model theoretic representation of unitary finite-rank perturbations. A short exposition on related Krein spaces and reproducing kernel Hilbert spaces is provided.
In Section \ref{s-STUnit} relationships between the family of spectra of the perturbation problem and the characteristic function are presented.

In the appendix Section \ref{s-APP} we take a moment to convey just the ideas behind several other well-deserving topics in the field. We refer to other literature for more information.

\subsection{Notation}\label{ss-notation}

We use different notation to help the reader distinguish between the unitary the self-adjoint setting.

In the self-adjoint setting, a rank-one perturbation of an operator $A$ will be denoted as $A_\gamma $, where $\gamma\in \R$. We will use ``boldface" $\bA\ciG $ for a finite-rank perturbation that is given by a self-adjoint matrix $(d\times d)$-matrix $\Gamma$. The real spectral measures for these cases will be referred to as $\mu_\gamma $ and $\bmu\ciG $ respectively. An additional superscript will be added when the trace of the matrix-valued spectral measures is required: $\bmu\ciG ^{\tr}$.
Also, the subscript will be entirely dropped when referring to objects corresponding to the unperturbed operator $\bA$, e.g.~$\bmu = \bmu\ci\OZ$, $F = F_0$, ${\bf F} = {\bf F_0}$, etc.

In the unitary setting, a rank-one perturbation of an operator $U$ will be denoted as $U_{\al}$, where $\alpha\in \T$. A finite-rank perturbation will be given by $\bU_{\al}$ with unitary $(d\times d)$-matrix $\alpha$. Notation similar to the self-adjoint setting will be used for the spectral measures, e.g.~$\mu_\al$ and $\bmu_\al$. Here, the subscript $\al$ indicates that we work with unitary perturbations. Characteristic functions and model spaces will be denoted in the rank-one case by $\te$ and $\cK_{\te}$, and in the finite-rank case by $\bte$ and $\cK_{\bte}$. 
Dropping the subscript again refers to objects that correspond to the unperturbed operator $\bU$, except this operator arises from using $\al = I$, e.g.~$\bmu = \bmu_I$, etc. We will simply write $I$ for the identity matrix, with the dimension inferred from context.

Spaces will be written in ``mathcal'' notation, e.g.~$\cH$ and $\cH_s(\bA)$. In particular, $\cD$ and $\cD_*$ refer to the deficiency spaces on the unitary side.

\section{Perturbation-theoretic background}\label{s-PRELIM}

We begin by presenting some central ideas from classical perturbation theory of self-adjoint operators, in order to better frame later discussions. 

A linear operator $A$ from a Banach space $\mathcal{X}$ to a Banach space $\mathcal{Y}$ is said to be \emph{compact} if the image $A(\mathcal{X}_1)$ of any bounded subset $\mathcal{X}_1\subset \mathcal{X}$ is relatively compact in $\mathcal{Y}$. Consider linear operators acting on a Hilbert space $\cH$. The class of compact operators $\mathscr{S}$ is then obtained by taking the closure of the set of finite-rank operators with respect to the operator norm topology. A characterization of the spectrum of self-adjoint operators that differ by a compact perturbation is available. Recall that the \emph{spectrum} of an operator $A$, denoted by $\sigma(A)$, is the closure of the set of all $\lambda \in \mathbb{C}$ for which operator $A-\lambda\OID$ is not invertible. The essential spectrum is the spectrum minus the isolated eigenvalues of finite (algebraic) multiplicity.

\begin{thm}[{von Neumann, see e.g.~\cite[Theorems 3 and 6 of Ch.~9]{BirmanSol-book_1987}}]\label{t-vN}
Let $A$ and $B$ be bounded self-adjoint operators. Then $B$ is compact if and only if the essential spectra of $A$ and $A+B$ are the same.
\end{thm}

In the self-adjoint setting, we can view compact operators as compact perturbations of the zero operator to see that compact operators are characterized as those whose only (possible) accumulation point of eigenvalues is the origin. A more refined standard definition restricts the speed at which the eigenvalues tend to $0$.
Namely, the \emph{von Neumann--Schatten classes}, $\mathscr{S}_p$, consist of compact operators whose sequence of singular values $\{s_k\}$ belongs to $ \ell^p$. Here, the singular values of an operator $T$ are defined as the eigenvalues of $|T| = (T^* T)^{1/2}$. Self-adjoint operators thus have the property that $s_k = |\lambda_k|$, where $\lambda$ is the sequence of eigenvalues.

\begin{thm}[Kato--Rosenblum, {\cite[Theorem 1]{Kato1957} and \cite[Theorem 1.6]{Rosenblum}}]\label{t-KR}
Let $A$ and $B$ be self-adjoint operators and assume $B\in \mathscr{S}_1$. Then the absolutely continuous parts of $A$ and $A+B$ are unitarily equivalent.
\end{thm}

Carey and Pincus \cite{CP} characterized trace class, $\mathscr{S}_1$, perturbations $A+B$ of $A$. Apart from leaving the absolutely continuous spectrum invariant, it must be possible to split the isolated eigenvalues of $A$ and $A+B$ as follows into three categories. The first and second categories are comprised of the eigenvalues of $A$ and of $A+B$, respectively, that have summable distance from the essential spectrum of $A$. The third category contains all remaining eigenvalues of $A$ and $A+B$. And there must exist a bijection $\f$ mapping those eigenvalues of $A$ in this category to those remaining eigenvalues of $A+B$ so that the sum of $|\lambda - \f(\lambda)|$ over all eigenvalues $\lambda$ of $A$ in this category is finite. In other words the remaining eigenvalues of $A$ have trace class distance to the remaining ones of $A+B$.

To emphasize a dichotomy, we mention that absolutely continuous spectrum can be destroyed by a Hilbert--Schmidt operator of arbitrarily small Hilbert--Schmidt norm:
\begin{thm}[Weyl--von Neumann, see e.g.~{\cite[p.~525]{katobook}}]
Let $A$ be a self-adjoint operator. 
For every $\eta>0$, there exists a self-adjoint operator $B$ with Hilbert-Schmidt norm less than $\eta$ so that $A+B$ has pure point spectrum.
\label{t-HSDestroy}
\end{thm}

Since the Hilbert--Schmidt norm dominates the standard operator norm, this means that the absolutely continuous spectrum may be unstable under arbitrarily small perturbations.

Theorem \ref{t-HSDestroy} was first proved by Weyl \cite{W} for compact perturbations and then for the smaller class of Hilbert--Schmidt perturbations by von Neumann \cite{1935}.
Extensions to normal operators and perturbations were proved by Berg \cite{1971}  for compact operators and by Voiculescu \cite{1979, 1981} for Hilbert--Schmidt perturbations. These results form the basis of $K$-homology theory, which studies the homology of the category consisting of locally compact Hausdorff spaces.

On the side, we mention Baranov \cite{spectral synthesis for normal operators} where a model representation and a spectral synthesis for rank-one perturbations of normal operators is achieved.

In order to avoid possible confusion, we spell out that we are not (at least not explicitly) reaching for a spectral synthesis, or other questions usually related to $K$-homology. Instead, we are primarily interested in spectral invariants and describing the spectral measure under perturbations.

\section{Aspects of self-adjoint rank-one perturbations}\label{s-rankone}

\subsection{Scales of Hilbert Spaces}\label{ss-scales}

When considering perturbations like Equation \eqref{e-rankone}, it is sometimes convenient to loosen our restrictions on the perturbation vector $\f$ to expand our possible applications, e.g.~to changing boundary conditions of differential operators. We say that the perturbation is \emph{bounded} when the vector $\f$ is from the Hilbert space $\cH$. The previous sections have dealt exclusively with bounded perturbations. If $\f\notin\cH$, we say the perturbation is \emph{singular}. These perturbations are significantly more complicated; it is imperative to ensure that the perturbation is well-defined in order to extend the tools that are presented  in Subsection \ref{ss-STR1P}. The description here roughly follows that of \cite{AK}.

Let $A$ be a self-adjoint (possibly unbounded) operator on a separable Hilbert space $\cH$. Consider the non-negative operator $|A|=(A^*A)^{1/2}$, whose domain coincides with the domain of $A$. Alternatively, if $A$ is bounded from below, the shifted operator $A+kI$, $k\in\RR$ sufficiently large, will provide a non-negative operator. We introduce a scale of Hilbert spaces.

\begin{defn}[\hspace{-1pt}{\cite[Section 1.2.2]{AK}}]\label{d-standardscale}
For $s\geq 0$, define the space $\cH_s(A)$ to consist of $\f$ from $\cH$ for which the $s$-norm
\begin{align}\label{d-scalenorm}
\|\f\|_s:=\|(|A|+I)^{s/2}\f\|_{\cH},
\end{align}
is bounded. 
The space $\cH_s(A)$ equipped with the norm $\|\cdot\|_s$ is complete. The adjoint spaces, formed by taking the linear bounded functionals on $\cH_s(A)$, are used to define these spaces for negative indices, i.e.~$\cH_{-s}(A):=\cH_s^*(A)$. The corresponding norm in the space $\cH_ {-s}(A)$ is thus defined by \eqref{d-scalenorm} as well. 
The collection of these $\cH_s(A)$ spaces will be called the \emph{scale of Hilbert spaces associated with the self-adjoint operator $A$}.
\end{defn}

It is not difficult to see that the spaces satisfy the nesting properties
\begin{align*}
\hdots\subset\cH_2(A)\subset\cH_1(A)\subset\cH=\cH_0(A)\subset\cH_{-1}(A)\subset\cH_{-2}(A)\subset\hdots,
\end{align*}
and that for every two $s,t$ with $s<t$, the space $\cH_t(A)$ is dense in $\cH_s(A)$ in the norm $\|\cdot\|_s$. Indeed, the operator $(A+1)^{t/2}$ defines an isometry from $\cH_s(A)$ to $\cH_{s-t}(A)$. In the rest of the subsection, we will use the brackets $\langle \fdot,\fdot\rangle$ to denote both the scalar product in the Hilbert space $\cH$ and the action of the functionals. For instance, if $\f\in\cH_{-s}(A)$, $\psi\in\cH_s(A)$, then 
\begin{align*}
\langle\f,\psi\rangle:=\big\langle(|A|+I)^{-s/2}\f,(|A|+I)^{s/2}\psi\big\rangle,
\end{align*}
where the brackets on the right hand side denote the scalar product.

Throughout the literature of other fields similar constructions occur under different names. For instance, the pairing of $\cH_1(A)$, $\cH$, and $\cH_{-1}(A)$ is sometimes referred to as a \emph{Gelfand triple} or \emph{rigged Hilbert space}. Also, when $A$ is the derivative operator, these scales are simply Sobolev spaces (with $p=2$). More details about Hilbert scales can be found in \cite{KP}. 

It is worth noting that these Hilbert scale are related to those generated by so-called left-definite theory \cite{LW02}. This theory employs powers of a semi-bounded self-adjoint differential operator to create a continuum of operators whereupon spectral properties can be studied. The theory can be applied to self-adjoint extensions of self-adjoint operators, which can be viewed as finite-rank perturbations, see e.g.~\cite{FFL, FL} and the references therein.

Rank-one perturbations of a given operator $A$ arise most commonly when the vectors $\f$ are bounded linear functionals on the domain of the operator $A$, so many applications are focused on $\cH_{-2}(A)$. Here, we only discuss the case $\f\in\cH_{-1}(A)$ for the sake of simplicity. However, references usually contain information on extensions to $\f\in\cH_{-2}(A)$, and information on the case when $\f\notin\cH_{-2}(A)$ can be found in \cite{DKS, Kurasov}.

\begin{rem}\label{r-H-n}
The case $\cH_{-1}(A)$ for the self-adjoint setting most closely aligns with unitary perturbations, see \cite[pp.~124--128]{LTSurvey}. It is not immediately clear how the more singular perturbations, $\cH_{-n}(A)$ for $n>1$, translate to the unitary side.
\end{rem}

\subsection{Spectral Theory of Rank-One Perturbations}\label{ss-STR1P}
A nice overview of what is now known as Aronszajn--Donoghue theory was given in \cite{SIMREV}. Extensions of Aronszajn--Donoghue theory to the case when the spectral measure is associated with a perturbation vector $\f\in\cH_{-2}(A)$ can found in \cite{AKK} and \cite{Kosh}, but here we take $\f\in\cH_{-1}(A)$ unless otherwise mentioned.
The results compare the spectral measures $\mu$ and $\mu_\gamma$ of the unperturbed and the perturbed operators and are expressed through the scalar-valued \emph{Borel} transform
\begin{align}\label{e-borel}
F_\gamma(z):=\int\ci\R \frac{d\mu_\gamma(t)}{t-z}\qquad\text{for }z\in \C\backslash\R,
\end{align}
which is abbreviated $F$ for $\gamma=0$. 
 
One of the standard identities at the heart of the theory is often referred to as the Aronszajn--Krein formula
$
F_\gamma(z)
=
F(z)/(1+\gamma F(z))$.
The distinction of whether or not a point has mass is encrypted in the functions $F$ and $\displaystyle G(x):= \int\frac{d\mu(t)}{(x-t)^2}$.

\begin{thm}[Aronszajn--Donoghue theory, e.g.~{\cite[Theorem 12.2]{SIMREV}}]\label{t-AD}
When $\gamma\neq 0$, the sets
\begin{align*}
S_\gamma
&=
\Bigl\{
x\in \R\Big|\lim_{y\to 0}F(x+iy) = -1/\gamma; G(x) = \infty
\Bigr\},\\
P_\gamma
&=
\Bigl\{
x\in \R\Big|\lim_{y\to 0}F(x+iy) = -1/\gamma; G(x) < \infty
\Bigr\}, \text{ and }\\
C
&=
\Bigl\{
x\in \R\Big|\lim_{y\to 0} \im F(x+iy) \neq 0
\Bigr\},
\end{align*}
contain spectral information of the perturbed operator $A_\gamma$ as follows:
\begin{itemize}
\item[(i)] For fixed $\gamma\neq 0$, the sets $S_\gamma$, $P_\gamma$ and $C$ are mutually disjoint.
\item[(ii)] Set $P_\gamma$ is the set of eigenvalues, and set $C$ ($S_\gamma$) is a carrier for the absolutely (singular) continuous measure, respectively.
\item[(iii)] For $\gamma\neq \beta$ the singular parts of $A_\gamma$ and $A_\beta$ are mutually singular.
\end{itemize}
\end{thm}

\begin{rem}
Set $X$ being a carrier for a measure $\tau$ means that $\tau(\R\backslash X) = 0$. Any (measurable) set that contains the support of a measure is also a carrier. Since we do not require a carrier to be closed, there may be carrier sets that are strictly contained in the support of a measure.
\end{rem}

The density function of the absolutely continuous measure and the pure point masses of $A_\gamma$ are completely described by the following result.

\begin{prop} Assume that $\gamma\neq 0$.
\begin{itemize}
\item[(i)] For $\lambda\in P_\gamma$ we have $\mu_\gamma(\{\lambda\}) = \frac{1}{\gamma^2 G(\lambda)}$.
\item[(ii)] The density function of the absolutely continuous part of $A_\gamma$ is given by
\[
\frac{d\mu_\gamma(x)}{dx}
=
\frac{1}{\pi}
\lim_{y\to 0^+}
\frac{\im F(x+iy)}{|1+\gamma F(x+iy)|^2},
\]
with respect to Lebesgue a.e.~$x\in \R$.
\end{itemize}
\end{prop}
We mention that the limit in part (ii) of the proposition exists with respect to Lebesgue a.e.~$x$. Indeed, by the Aronszajn--Krein formula $\frac{\im F}{|1+\gamma F|^2} = \im F_\gamma$, and $F_\gamma$ is analytic on the upper half-plane.

A  characterization of the singular continuous part of  $A_\gamma$ has been sought after but is still outstanding. Only partial results have been established. Instead of elaborating on the details here, we refer the reader to \cite{cimaross, LT09, SIMREV} and the references therein. We also point the reader to \cite{LTSurvey} for a discussion of, and references for, the question: ``How unstable can the singular spectrum become?''

The measures $\mu$ and $\mu_\gamma $, which are the spectral measures associated with rank-one perturbations of self-adjoint operators, are associated with scalar Nevanlinna--Herglotz functions. These functions are analytic self-maps of the upper half plane $\CC_+$ and possess the Nevanlinna--Riesz--Herglotz representation 
\begin{align*}
    \widetilde F(z)=c+dz+\int_{\RR}\left(\dfrac{1}{t-z}-\dfrac{t}{1+t^2}\right) d\mu(t),
\end{align*}
and $\mu$ is a measures which satisfies the decay condition $\int_{\RR}(1+t^2)^{-1}d\mu(t)<\infty$. The examples that use $\widetilde F$ are more singular $\cH_{-2}(A)$ perturbations.
In order to give the reader additional intuition about these measures, we include some examples from \cite[App.~A]{GT}.

\vspace{1em}
\begin{center}
\renewcommand*{\arraystretch}{1.5}
 \begin{tabular}{||c|c|c||} 
 \hline
 \cite[App.~A]{GT} & Borel transform & Spectral Measure $d\mu(t)$\\
 \hline\hline
  Eq.~(5) & $F(z)= -1/z$ & $\delta\ci{\{0\}}(t)dt$ \\
 \hline
 Eq.~(6) & $\widetilde F(z)= \ln(z)$ & $\chi\ci{(-\infty,0)}(t)dt$  \\
 \hline
 Eq.~(7) & $\widetilde F(z)= \ln(-1/z)$ & $\chi\ci{(0,\infty)}(t)dt$ \\
 \hline
Eq.~(8) & $\widetilde F(z)= z^r-\cos(\frac{r\pi}{2})$ & $|t|^r\pi^{-1}\sin(r\pi)\chi\ci{(-\infty,0)}dt$, $r\in (0,1)$ \\
 \hline
 Eq.~(10) & $\widetilde F(z)= \tan(z)$ & $\sum_{n\in\ZZ}\delta\ci{\{n\pi\}}(t)dt$  \\
 \hline
 Eq.~(17) & $F(z)= \ln\left(\frac{z-t_1}{z-t_2}\right)$ & $\chi\ci{[t_1,t_2]}(t)dt$ with $t_1<t_2$  \\ [1ex] 
 \hline
\end{tabular}\\\vspace{.2cm}
\small{Examples of Nevanlinna--Herglotz functions and corresponding spectral measures. The examples Eq.~(6)--(8) use the principal value of the logarithm. The integration variable is $\lambda$. The first column contains references to equations in  \cite[App.~A]{GT}. Other examples and their sources can also be found there.}
\end{center}

\section{Aspects of Unitary Rank-One Perturbations and Model Theory}\label{s-introunit}

Consider the unitary rank-one perturbation problem given by equation \eqref{d-sUalpha}. Let $\mu_\al$ be the spectral measure of $U_\al$ with respect to the $*$-cyclic vector $\f$, which is simultaneously also $*$-cyclic for $U_\al$ for all $\al\in \T$. Then, the Spectral Theorem says $U_\al$ can be represented by the operator that acts via multiplication by the independent variable on the space $L^2(\mu_\al)$. 

The operator $U_0$ is well-known to be a completely non-unitary contraction, i.e.~it is not unitary on any of its invariant subspaces. Therefore, it (and hence the family of measures $\{\mu_\al\}$) corresponds to the compression of the shift operator in a model representation associated with a characteristic function $\te$. Studying the intricacies of these model representations emerges as one of the main strategies in this field.

Model spaces are subspaces of a weighted $L^2$-space, of which we discuss several: the one by Clark, which resembles a simplified Sz.-Nagy--Foia\c s model; the one by de Branges--Rovnyak which was e.g.~studied by the Sarason school; and an overarching description of model theory developed by Nikolski--Vasyunin. This final formulation essentially incorporates the former ones by choosing an appropriate weight function.

\subsection{Aleksandrov--Clark Theory and Sz.-Nagy--Foia\c s Model for Perturbations with Purely Singular Spectrum}\label{ss-SzNF}

A seminal paper by Clark \cite{Clark} laid the foundation that connects rank-one perturbations with reproducing kernel Hilbert spaces. The field has since grown into what is now known as Aleksandrov--Clark theory, honoring the deep insights gained by Aleksandrov about Clark measures -- especially in the presence of an absolutely continuous component. A nice exposition of Aleksandrov--Clark theory can be found in \cite{cimaross}, which we mostly follow along with in this section.
We refer readers interested in a more general exposition of the Sz.-Nagy--Foia\c s model spaces to \cite{SzNF2010}. For roughly the second half of this subsection, we work with characteristic functions that are inner, or equivalently, within the Clark setting of purely singular spectral measures.

For an analytic function $\te:\DD\to\DD$ and a point $\al\in\TT$, the function 
\begin{align}\label{e-ualpha}
    u_\al(z):=\Re\left(\dfrac{\al+\te(z)}{\al-\te(z)}\right)=\dfrac{1-|\te(z)|^2}{|\al-\te(z)|^2},
\end{align}
is positive and harmonic on $\DD$. For each $\al$, a theorem by Herglotz \cite{H} says this function corresponds uniquely to a positive measure $\mu_\al$ with $u_\al=P\mu_\al$. Here, $P\mu_\al=\int_{\TT}\frac{1-|z|^2}{|\zeta-z|^2}d\mu_{\al}(\zeta)$ is the Poisson integral of $\mu_\al$. 

We let
$\cA_{\te}:=\{\mu_{\al}:\al\in\TT\}$ 
denote the family of measures associated with the function $\te$. We will call $\cA_{\te}$ the family of {\em Clark measures} of $\te$ when $\te$ is an inner function, i.e.~a bounded analytic function with unit modulus a.e.~on $\TT$. Note that when $\te$ is a general analytic self-map of the disk, the family $\cA_{\te}$ is usually referred to as the {\em Aleksandrov--Clark measures} of $\te$.

With the Herglotz transformation $(H\mu)(z)=\int_{\TT}\frac{\zeta+z}{\zeta-z}d\mu(\zeta)
$ of a measure $\mu = \mu_1$, it can easily be verified that the function
\begin{align}\label{e-HerglotzRep}
    \te(z):=\dfrac{(H\mu)(z)-1}{(H\mu)(z)+1},
\end{align}
is an analytic self map of the disk. The condition $\te(0)=0$ is equivalent to each $\mu_{\al}\in\cA_{\te}$ being a probability measure {\cite[Proposition 9.1.8]{cimaross}}. 

These Clark measures can be used to describe the unitary perturbations of an important operator. To do so, we define the shift operator $S:H^2\to H^2$ by $(Sf)(z)=zf(z)$, where $H^2 = H^2(\D)$ denotes the Hardy space. Likewise, for later, we define the backward shift operator to be $(S^*f)(z)=\frac{f(z)-f(0)}{z}$. Beurling's Theorem \cite{Beurling} then says that the $S$-invariant subspaces of $H^2$ are exactly those that can be written as $\te H^2$ for some inner function $\te$.

In order to take advantage of this relationship, we now assume $\te$ is an inner function with $\te(0)=0$. Assuming that $\te(0)=0$ is not essential, but rather a convenience. Sometimes we will refer to such functions as {\em characteristic} functions. Given such a $\te$, the \emph{Sz.-Nagy-Foia{\c s} model space} \cite{SzNF2010} can then be defined as
\begin{align}\label{d-Kte}
    \cK_{\te}:=H^2\ominus \te H^2.
\end{align}
Beurling's Theorem further implies that $S^*$-invariant subspaces of $H^2$ are simply model spaces $\cK_{\te}$ corresponding to some inner $\te$.

On the side, we mention two major advances in complex analysis:\\
(i) Douglas--Shapiro--Shields \cite{DSS} have shown that for $f\in H^2$, $f\in \cK_{\te}$ if and only if the meromorphic function $f/\te$ on $\DD$ has a pseudo-continuation to a function $\widetilde{f_{\te}}\in H^2(\C\setminus\overline{\D})$ with $\widetilde{f_{\te}}(\infty)=0$ (also see \cite[Theorem 8.2.5]{cimaross}). The analogous result was also shown there to hold for conjugate pairs of $H^p$ spaces.\\
(ii) A milestone has been achieved with the Ahern--Clark Theorems \cite{AhernClark1, AhernClark2} with respect to understanding when the Julia--Carath\'eodory angular derivative exists. This result was generalized by Fricain--Mashreghi \cite{FricMash} to a characterization based on the existence of radial limits for higher derivatives. Also see the survey by Garcia--Ross \cite[Theorem 6.11]{GarciaRoss} for a summary.

Moving on with our program, let $P_{\te}$ be the orthogonal projection of $H^2$ onto $\cK_{\te}$. The \emph{compression of the shift operator} is thus defined as
\begin{align*}
    S_{\te}=P_{\te}S|\ci{\cK_{\te}}.
\end{align*}
This allows us to write the family of rank-one perturbations on $\cK_\te$:
\begin{align}
V_{\al}f=S_{\te}f+\al\left\langle f,\frac{\te}{z}\right\rangle \ID, \qquad\text{with}\qquad\al\in\TT.
\end{align}
In particular, the following theorem of Clark says that these are the only unitary rank-one perturbations of $S_{\te}$.

\begin{thm}[Clark {\cite[Remark 2.3]{Clark}}]
Any operator $X$ that is both unitary and a rank-one perturbation of $S_{\te}$ can be written as $X=V_{\al}$ for some $\al\in\TT$.
\end{thm}

Let $\mu_{\al}$ be the Clark measure associated with the inner function $\te$ and the point $\al\in\TT$. Since $V_{\al}$ is a cyclic unitary operator, the spectral theorem says that $V_{\al}$ can be represented as multiplication by the independent variable on some $L^2(\nu)$ space. It turns out that the space $L^2(\nu)$ can be canonically identified with $L^2(\mu_{\al})$. Let $M$ be the operator on the space $L^2(\mu_{\al})$ acting via multiplication by the independent variable. 
Then, the unitary operator that intertwines, $C_{\al}M=V_{\al}C_{\al}$, and maps the constant function $\ID\in L^2(\mu_\al)$ to some vector in the defect space $\Ran{(I-S_\te^*S_\te)^{1/2}}$ is called the \emph{adjoint Clark operator}. It is given by the
\emph{normalized Cauchy transform}
\begin{align*}
    C_{\al}:
    L^2(\mu_\al)\to \text{Hol}(\DD)
    \qquad\text{with}\qquad
    (C_{\al}g)(z):=\dfrac{K(gd\mu_\al)}{K\mu_\al},
\end{align*}
where $K$ is the Cauchy transform
    $(K\nu)(z)=\int_{\TT}\frac{d\nu(\zeta)}{1-z \bar\zeta}.$
    The Clark operator is often denoted by $\Phi$ in literature, so that $C_\al = \Phi^*$.

These representations gives us access to spectral information regarding the Clark family $\{\mu_\alpha\}$, $\alpha\in \T$. 

\begin{thm}[see e.g.~{\cite[Proposition 9.1.14]{cimaross} and \cite[Proposition 8.3]{GarciaRoss}}]\label{t-STR1Unit}
In the above setting we have:
\begin{enumerate}
\item $(d\mu_\al)\ti{ac} = u_\al dm$ (with $u_\al(z) = (1-|\te(z)|^2)|\al-\te(z)|^{-2}$ as in  \eqref{e-ualpha}).
\item $\mu_\alpha\perp\mu_\beta$ for all $\alpha\neq\beta$, $\beta\in\T$.
\item $\mu_\alpha$ has a point mass at $\zeta\in \T$ if and only if $\te(\zeta)=\alpha$ and $|\te'(\zeta)|<\infty$. In that case this point mass is given by $\mu_\alpha(\{\zeta\}) = |\te'(\zeta)|^{-1}$.
\item The set $\{\zeta\in \T:\lim_{r\to 1^-}\te(r\zeta) = \alpha\}$ is a carrier for $\mu_\alpha$. (Recall that $\mu_\alpha$ is purely singular in the Clark setting.)
\end{enumerate}
\end{thm} 

This result is in direct correspondence with the Aronszajn--Donoghue Theorem \ref{t-AD} above. Also, observe that a point mass equals the reciprocal magnitude of the derivative of the Borel transform in the self-adjoint setting, and of the Cauchy transform in the unitary setting. In fact, \cite[Item (1) of Corollary 9.1.24]{cimaross} offers a finer carrier of the singular spectrum in terms of the lower Dini derivative of $\mu_\al.$

The Sz.-Nagy--Foia\c s  representation simplifies to the setting described in this subsection precisely when operator $V_1$ has no absolutely continuous part (or, equivalently, when the characteristic function $\te$ is inner). This poses a significant restriction. The de Branges--Rovnyak model is an alternative representation of the situation under weaker conditions. In the most general Aleksandrov--Clark situation, one is required to deal with the full two-storied 
Sz.-Nagy--Foia\c s model space
\[
 \cK_\te=
 \left( \begin{array}{c} H^2 \\ \clos\Delta L^2 \end{array}\right) 
 \ominus 
  \left( \begin{array}{c} \theta \\ \Delta \end{array}\right) H^2,
\]
instead of just the first component as in \eqref{d-Kte}. The defect function $\Delta$ is $\Delta(z) = (1-\theta(z)^*\theta(z))^{1/2}$ for $z\in \T$. For further reference, see \cite[Section 1.3.5]{Nikolski}.

\subsection{de Branges--Rovnyak Model and Perturbations in the Extreme Case}\label{ss-deBR}

In this subsection, we assume that the characteristic function $\te$ is an {\em extreme} point, i.e.~ that $\int\ci\T \ln(1-|\te(z)|)dm(z)=-\infty$. It is well-known that $\te$ extreme if and only  $L^2(\mu) = H^2(\mu)$ for the corresponding Aleksandrov--Clark measure $\mu = \mu_1$. This situation is ideal for the de Branges--Rovnyak model space, as it now reduces from two components
\[
\cK_\theta = \left\{
\left(\begin{array}{c} g_+ \\ g_-\end{array} \right) \,:\ g_+\in H^2
,\  g_- \in H^2_-
,\  g_- - \theta^* g_+ \in \Delta L^2
\right\}  
\]
to a one component space. Here we used the notation $H_-^2:=L^2\ominus H^2.$

We describe the reduced one-component de Branges--Rovnyak model space: So, assume $\te:H^2
\to H^2
$ is extreme. Then the de Branges--Rovnyak model space $\cH(\te)\subset H^2
$ consists of functions in the range space of the defect operator, i.e.~$\cH(\te) = (\OID - |\te|^2)^{1/2}H^2
.$
The canonical norm on this space is the \emph{range norm} which arises by taking the minimal norm of the pre-image of an element from $\cH(\te)$. Much of the success of this approach is based upon the fact that $\cH(\te)$ is a reproducing kernel Hilbert space with reproducing kernel $k^\te_w(z) = \frac{1-\overline{\te(w)}\te(z)}{1-\bar w z}.$ The deep structure of this space is the focus of \cite{SAR}. Here we only mention a few items relevant to perturbation theory. We will omit other interesting topics such as multipliers of $\cH(\te)$, the theory regarding the Julia--Carath\'eodory angular derivatives and Denjoy--Wolff points --- all of which are detailed in \cite{SAR}. 

The connection with the corresponding Aleksandrov--Clark measure $\mu$ is made through equation \eqref{e-HerglotzRep}, see e.g.~\cite[Chapter III]{SAR}. Much of the development in this area is attributed to the dissertation of Ball \cite{BallDiss}. For instance, it was shown there that the measure $\mu$ has an atom at a point $z_0\in\T$ if and only if the function $\frac{\te(z)-1}{z-z_0}$ belongs to $\cH(\te)$, see e.g.~\cite[Section (III-12)]{SAR}.

\subsection{General perturbations and Nikolski--Vasyunin Model Theory}\label{ss-NV}

Not all rank-one perturbations satisfy any of the conditions under which we can use the representations detailed in Subsections \ref{ss-SzNF} and \ref{ss-deBR}. Model theory for unitary perturbations in the general setting is much more complicated. Instead of a one-story model space, the general setting requires a two-story model space. While this description is superior in abstraction and admits more general settings, the models discussed in Subsections \ref{ss-SzNF} and \ref{ss-deBR} have provided many deep insights over the years.

An overarching treatise of the Sz.-Nagy--Foia\c s, the de Branges--Rovnyak model space and other model spaces (e.g.~the one studied by Pavlov) was achieved by Nikolski--Vasyunin \cite{Nikolski, Nik-Vas_model_MSRI_1998, Nik-Vas_TwoModels_1985, Nik-Vas_FunctModels_1989}. There, a general so-called \emph{transcription free} model space was introduced as a subspace of a (possibly) two-storied weighted space $L^2(\fD_*\oplus \fD, W)$ on the unit circle. Here, the defect spaces of contraction $V_0$ are given by $\fD =  \clos\Ran{(I - S_\te^* S_\te)^{1/2}}$ and $\fD_*= \clos\Ran{(I - S_\te S_\te^*)^{1/2}}$.  We also note that the defect spaces $\cD$ and $\cD_*$ were identified with $\T$ in Subsections \ref{ss-SzNF} and \ref{ss-deBR}. This $L^2$ space then reduces to the Sz.-Nagy--Foia\c s, the de Branges--Rovnyak, the Pavlov model spaces, and other transcriptions by making specific choices of the weight $W$.
The connection to rank-one perturbations comes from the dependence of this $W$ on the characteristic function $\te$.

General rank-one perturbations were studied in Liaw--Treil \cite{LT17-Clark}. This subject is included in the lecture notes by Liaw--Treil \cite{LTSurvey} on the relationship between rank-one perturbations and singular integral operators. Instead of repeating large chunks of information here, we refer the reader to those lecture notes.

\section{Self-Adjoint Finite-Rank Perturbations}\label{s-SAFRP}

We adapt the self-adjoint finite-rank setup given in \eqref{d-bAGamma} to account for singular perturbation vectors, which are useful in many applications. We mostly follow along with \cite[Ch.~3]{AK}. We begin by defining the coordinate operator $\bB:\C^d\to \cH_{-2}(\bA)$ that takes the standard basis $\{{\bf e}_k\}_{k=1}^d\subset\C^d$ to $\{\f_k\}_{k=1}^d\subset\cH_{-2}(\bA)$. Note that we are changing notation slightly from Section \ref{s-Intro}, as we used to think of $\bB$ as an operator $\bB:\C^d\to \Ran{\bB}$ that was invertible. As before, we assume without loss of generality the invertibility of $\bB$ on its range.

Consider finite-rank perturbations of a self-adjoint operator $\bA$ on the separable Hilbert space $\cH$ given by
\begin{align}\label{e-newSAFR}
\bA\ciG =\bA+\bB\Gamma\bB^*,
\end{align}
where $\Gamma$ is a Hermitian $d\times d$ matrix and the operator $\bB\Gamma\bB^*$ is an operator of rank $d$ from the Hilbert space $\cH_{2}(\bA)$ to the Hilbert space $\cH_{-2}(\bA)$. Note that we can assume without loss of generality that the matrix $\Gamma$ is invertible. If $\Gamma$ is not invertible, then the orthogonal complement to the kernel of the operator $\Gamma$ yields a finite-rank operator of rank strictly less than $d$ determined by a non-degenerate Hermitian matrix. 

The vectors $\f_k$ can be thought of as modifying the domain of $\bA$ by $d$ dimensions that are in $\cH_{-2}(\bA)$. However, to ensure that each of these vectors are non-degenerate and adding new dimensions, we will call the set of vectors $\f_k\in\cH_{-2}(\bA)\backslash\cH$, $k=1,\dots,d$, {\em $\cH$-independent} if and only if the equality
\begin{align*}
    \sum_{k=1}^d c_k\f_k\in\cH, \text{ } c_k\in\CC,
\end{align*}
implies $c_1=c_2=\cdots=c_d=0$. If a desired set is not $\cH$-independent, then the matrix $\bB\Gamma\bB^*$ will not be invertible and define a degenerate perturbation of rank strictly less than $d$. For this reason, we consider only $\cH$-independent perturbations.

\subsection{Singular Finite-Rank Perturbations}

The operator $\bA\ciG $ on the domain $\Dom(\bA)$ is symmetric as an operator acting from $\cH_2(\bA)=\Dom(\bA)$ to $\cH_{-2}(\bA)$. The self-adjoint operator given by equation \eqref{e-newSAFR} coincides with one of the self-adjoint extensions of the operator $\bA^0$ equal to the operator $\bA$ restricted to the domain
\begin{align*}
    \Dom(\bA^0)=\Dom(\bA)\cap\Ker(\bB\Gamma\bB^*).
\end{align*}
On the side we mention that $\Ker(\bB\Gamma\bB^*)=\Ker(\bB^*),$ because we are assuming $\Gamma$ to be invertible and $\cH$-independence of $\f_k$.

\begin{lem}[{\cite[Lemma 3.1.1]{AK}}]\label{l-newminimal}
Suppose that the vectors $\f_k\in\cH_{-2}(\bA)\backslash\cH$, $k=1,\dots,d$, are $\cH$-independent and form an orthonormal system in $\cH_{-2}(\bA)$. Then the restriction $\bA^0$ of the operator $\bA$ to the domain $\Dom(\bA^0)$ is a densely defined symmetric operator with the deficiency indices $(d,d)$.
\end{lem}

Note that the vectors $\f_k$ having unit norm in $\cH_{-2}({\bf A})$ is not a restriction, as every $\cH$-independent system $\{\f_k\}$ can be orthonormalized. We assume unit norm in the following discussions and results.

If we let the vectors $\f_k$, $k=1,\dots,d$ be $\cH$-independent, all vectors $\psi\in\Dom({\bf A}^{0*})$ can be represented as:
\begin{align}\label{e-newgreens}
    \psi=\widehat{\psi}+\sum_{k=1}^d \left( a_{+k}(\psi)(\bA-iI)^{-1}\f_k+a_{-k}(\psi)(\bA+iI)^{-1}\f_k\right),
\end{align}
where $\widehat{\psi}\in\Dom(\bA^0)$, $a_{\pm}(\psi)\in\CC$. 

The theory of self-adjoint extensions of symmetric differential operators, commonly referred to as Glazman--Krein--Naimark theory \cite{AG, N}, should be compared to this setup. The $\Dom(\bA^0)$ should be thought of as a ``minimal'' domain for the operator $\bA$, as the domain is unaffected by the perturbation $\bB\Gamma\bB^*$ and will be contained in the domains of all extensions. Likewise, the ``maximal'' domain is represented by $\Dom(\bA^{0*})$ and equation \eqref{e-newgreens} is a modified version of the classical von Neumann's formula (the maximal domain is the direct sum of the minimal domain and the defect spaces). The key space $\Dom(\bA^{0*})$ should thus be considered as a finite dimensional extension of the space $\cH_2(\bA)$ in the sense that $\Dom(\bA^{0*})$ is isomorphic to the direct sum of $\cH_2(\bA)$ and $\CC^d$.

We also emphasize that the spaces $\bA^{0}$, defined via Lemma \ref{l-newminimal}, and $\bA^{0*}$, are dependent on the choice of the vectors $\{\f_k\}_{k=1}^d$. We can thus formulate a second scale of Hilbert spaces
\begin{align*}
    \Dom({\bf A})=\cH_2({\bf A}) \subset \Dom(\bA^{0*}) \subset \cH \subset \Dom(\bA^{0*})^* \subset \cH_{-2}({\bf A})=\Dom({\bf A})^*,
\end{align*}
which is constructed using both the operators ${\bf A}$ and $\bB\Gamma\bB^*$. The norms in $\cH_{-2}(\bA)$ and $\cH_2(\bA)$ are the standard norms from Definition \ref{d-standardscale}. We avoid most of the specific properties of these spaces and operators, but point out that the norm in the space $\Dom(\bA^{0*})^*$ is listed in {\cite[Equation (3.11)]{AK}}, near other pertinent facts.

\subsection{Self-Adjoint Extensions}\label{ss-saextensions}

The self-adjoint finite-rank perturbation given by \eqref{e-newSAFR} can be adapted as an application to self-adjoint extension theory. Namely, self-adjoint extensions of the operator $\bA^0$ are parametrized by $d\times d$ unitary matrices by the classical Glazman--Krein--Naimark theory \cite{AG, N}. Let $V$ be such a matrix and the vector notation $\vec{a}_{\pm}\equiv\{a_{\pm}\}_{k=1}^d$ denote the coefficients from equation \eqref{e-newgreens}. The corresponding self-adjoint operator $\bA(V)$ coincides with the restriction of the operator $\bA^{0*}$ to the domain 
\begin{align}\label{e-saematrix}
    \Dom(\bA(V))=\{ \psi\in\Dom(\bA^{0*}) ~:~ -V\vec{a}_-(\psi)=\vec{a}_+(\psi)\}.
\end{align}
We present an explicit connection between $V$ and $\Gamma$ in Lemma \ref{l-VGamma} below.

The extension given by the matrix $V=I$ coincides with the original operator $\bA$. This case is handled by classical self-adjoint extension theory. However, when the perturbing vectors $\{\f_k\}_{k=1}^d$ belong to $\cH_{-1}(\bA)$, descriptions of the corresponding domains become more difficult.

\begin{thm}[{\cite[Theorem 3.1.1]{AK}}]
Let $\f_k\in\cH_{-1}(\bA)\backslash\cH$ be an $\cH$-independent basis such that $\langle(\bA-iI)^{-1}\f_j,(\bA+iI)^{-1}\f_k\rangle=\delta_{jk}$, and let $\Gamma$ be a Hermitian invertible matrix. Then the self-adjoint operator $\bA\ciG =\bA+\bB\Gamma\bB^*$ is the self-adjoint restriction of the operator $\bA^{0*}$ to the following domain
\begin{align*}
    &\Dom(\bA\ciG )\\
    =&\{\psi\in\Dom(\bA^{0*}) ~:~\vec{a}_+(\psi)=-(\Gamma^{-1}+{\bf F}(i))^{-1}(\Gamma^{-1}-{\bf F}^*(i))\vec{a}_-(\psi)\},
    \end{align*}
where ${\bf F}(i)=\bB (\bA\ciG-i I)^{-1}\bB^*$.
\end{thm}

The notation ${\bf F}(i)$  comes from the Borel transform, which we focus on in Section \ref{s-SAST}.

We have ${\bf A}_0={\bf A}$ when $\Gamma=0$. Further note that the matrix $V=(\Gamma^{-1}+{\bf F}(i))^{-1}(\Gamma^{-1}-{\bf F}^*(i))$ is unitary. Hence, the theorem says that if the vectors $\f_j$ and the desired perturbation $\Gamma$ are known, then the domain of the self-adjoint extension can be written via the explicit unitary matrix $V$, as in the classical theory. 

However, this leads to the natural question: Given the domain of a self-adjoint extension in terms of $V$, can we recover the perturbation $\Gamma$ responsible for this domain? The answer is given by the following result.

\begin{lem}[{\cite[Lemma 3.1.2]{AK}}]\label{l-VGamma}
Let $\f_k\in\cH_{-1}(\bA)\backslash\cH$, $k=1,\dots,d$, be an $\cH$-independent orthogonal system. If $$\det
\left( V+[iI+\RE({\bf F}(i))]^{-1}[iI-\RE({\bf F}(i))]
\right)
\neq 0$$ then the operator $\bA^{0*}$ restricted to the domain of functions $$\{\psi\in\Dom(\bA^{0*}) ~:~ -V\vec{a}_-(\psi)=\vec{a}_+(\psi)\}$$ is a finite dimensional additive perturbation of the operator $\bA$. In particular, the Hermitian invertible matrix $\Gamma$ is given by
\begin{align*}
    \Gamma=\left(-\RE({\bf F}(i))+i(I-V)^{-1}(I+V)\right)^{-1}.
\end{align*}
\end{lem}

The last formula necessitates the analysis of whether $I-V$ is invertible. This distinction is handled in the proof, where it is determined that if $I-V$ is not invertible, then there is a degeneracy in the choice of the vectors $\f_k$. This means that the set of vectors $\{\f_k\}_{k=1}^d$ contains extra elements because we can find a new set of elements $\{\f_k^*\}_{k=1}^{d^*}$, $d^*<d$, such that the corresponding matrix $V^*$ has a trivial eigensubspace. 

The description of domains of self-adjoint extensions resulting from finite-rank perturbations with vectors from $\cH_{-2}(\bA)$ are much more involved (see e.g.~\cite{AK2}), and while very interesting in their own right, fall outside the scope of our discussion.

\subsection{Some Applications of Singular Finite Rank Perturbations}

The singular finite-rank perturbation setup employed in this section has a wide array of applications. Perhaps, their most common uses include point interactions for differential operators via connections to distribution theory and singular potentials of Schr\"odinger operators. This is immediately evident from the rank-one case when considering changing boundary conditions of regular Sturm--Liouville operators, see \cite[Section 11.6]{Simon}. 

Several contributions to the finite-rank case can be found in \cite{AK}. These include the analysis of operators with generalized delta interactions to achieve both spectral and scattering results. It is also possible to consider infinite-rank perturbations, under some simplifying assumptions, to help approach problems given by two-body, three-body and few-body models.

Finally, we should mention that singular perturbations can be transcribed into the theory of rigged Hilbert spaces, i.e.~\cite{KD}. This theory places a larger emphasis on properties of singular quadratic forms, which can also describe self-adjoint extensions. Specific extensions, such as the Friedrichs or von Neumann--Krein cases, are sometimes easier to formulate in this context. Various aspects of spectral theory for singular finite-rank perturbations of self-adjoint operators are detailed in \cite[Section 9]{KD}.

\section{Spectral Theory of Self-Adjoint Finite-Rank Perturbations}\label{s-SAST}

Consider the family of finite-rank perturbations $\bA\ciG =
\bA
+\bB\Gamma\bB^*$, see \eqref{d-bAGamma}, with cyclic subspace $\Ran{\bB}$. It is well-known that $\Ran{\bB}$ is then also cyclic for $A\ciG$ for all symmetric $\Gamma$. For simplicity let us focus on bounded perturbations in this section. By the Spectral Theorem, this perturbation family corresponds to a family of matrix-valued spectral measures $\bmu\ciG$ through
\[
\bB^*(\bA\ciG - z I)^{-1}\bB
=
\int\ci{\R}\frac{d\bmu\ciG(t)}{t-z}
\qquad\text{for }z\in \C\setminus\R.
\]
The right hand side is the matrix-valued Borel transform, ${\bf F}\ciG(z): = \int\ci{\R}(t-z)^{-1}d\bmu\ciG(t).$ 
We obtain the \emph{scalar spectral measures} $\mu\ciG$ by taking the trace of $\bmu\ciG$. This trace is a scalar-valued measure which recovers the spectrum of $\bA\ciG$ via $\sigma(\bA\ciG) = \supp \mu\ciG.$ However, to access more subtle information, we formulate some of the results of the field we define the family of matrix-valued functions $W\ciG$ by $d\bmu\ciG(t) = W\ciG (t)d\mu\ciG(t)$. Finally, we arrive at $(W\ciG)\ti{ac}:= d\bmu\ciG / dx$ by taking a component-wise Randon--Nikodym derivative.

\subsection{Absolutely Continuous Spectrum and Scattering Theory}

The unitary equivalence of the absolutely continuous spectrum of operators that differ by a finite-rank perturbations is available through simply applying the Kato--Rosenblum Theorem \ref{t-KR}. In the more general setting of compact perturbations, the standard proof relies on the existence of the wave operators. Namely, let $P\ti{ac}$ denote the orthogonal projection from the Hilbert space onto the absolutely continuous part of $A$. For self-adjoint $A$ and compact self-adjoint $K$ it was shown that the strong operator topology limit of $e^{it(A+K)}e^{-itA} P\ti{ac}$ exists (see \cite[Theorem 1.6]{Rosenblum} and \cite[Theorem 1]{Kato1957}), which in turn yields the Kato--Rosenblum theorem.

For finite-rank perturbations, a quicker proof of unitary equivalence of the absolutely continuous spectrum is available. This proof uses the Aronszajn--Krein relation 
\begin{align}\label{e-AKMatrix}
{\bf F}\ci\Gamma = (\OID + {\bf F}\Gamma)^{-1} {\bf F}= {\bf F}(\OID + {\bf F}\Gamma)^{-1},
\end{align}
of the matrix-valued Borel transforms and was first discovered by Kuroda \cite{Kuroda1963}. Via efficient notation, and in a slightly different language, Liaw and Treil \cite[Appendix A.1]{LTJST} present this proof in a format appropriate for a graduate course.

Of course, scattering theory is able to give us more information by relating how wave operators, e.g.~ $s-\lim_{t\to\infty}e^{it(A+K)}e^{-itA} P\ti{ac}$, and their packets are affected by the perturbation. For an interesting exposition of scattering theory for finite-rank perturbations confer, e.g.~\cite[Ch.~4]{Kurasov}. Applications in Mathematical Physics can also be found in \cite[Ch.~5--7]{Kurasov}. Alternatively, the scattering theory of finite-rank perturbations can be analyzed using boundary triples, see e.g.~\cite{BMN}.

Validating the observation that the behavior of the absolutely continuous spectrum is one of the easier objects to capture, we conclude this subsection with its full perturbation theoretic characterization. The density of the matrix-valued spectral measure of the perturbed operator $(W\ciG)\ti{ac}$ is determined (see \cite[Lemma A.3]{LTJST}) in terms of that of the unperturbed operator $W\ti{ac}$ by
\[
	\left(W\ciG\right)\ti{ac}(x) =
	\lim_{y\to 0^+
	}(\OID +{\bf F}(x+iy)^* \Gamma )^{-1}
	W\ti{ac}(x)
	\lim_{y\to 0^+
	}(\OID +\Gamma {\bf F}(x+iy))^{-1},
	\]
	with respect to Lebesgue a.e.~$x\in \R$.

In Equation \eqref{e-ACinTheta} below, we also include a full description of the perturbed operator's matrix-density in terms of the matrix characteristic function of a corresponding model representation.

\subsection{Vector Mutually Singular Parts}

As evidenced by much research in the field, working with the singular spectrum will require a more subtle analysis than is necessary for the absolutely continuous part. From a naive perspective, the task at hand is to attempt to obtain some information about non-tangential boundary values $z\to \lambda$ of matrix-valued analytic functions on $\D$ for  $(\mu\ciG)\ti{s}$-a.e.~$\lambda\in \T$. As we discuss in Section \ref{ss-PoltThm}, Poltoratski's Theorem does not hold in the matrix-valued setting. Yet some positive results prevail.

Recall the Aronszajn--Donoghue Theorem, which states the mutual singularity of the singular parts under rank-one perturbations, see item (iii) of Theorem \ref{t-AD}. For finite-rank perturbations it is easy to construct examples for which two different perturbed operators have the same eigenvalue by taking direct sums of rank-one perturbations. The eigenvalues of the different components are completely independent from one another. Hence, a literal extension of this Aronszajn--Donoghue result cannot be true for the scalar-valued spectral measure. Through defining a vector-valued analog of the mutual singularity of matrix measures,  Liaw--Treil \cite[Theorem 6.2]{LTJST} achieved such a generalization of the Aronszajn--Donoghue Theorem. The scalar-valued spectral measures are also restricted:
\begin{thm}\cite[Theorem 6.3]{LTJST}
Fix a singular scalar Radon measure $\nu$, and $d\times d$-matrices $\Gamma>0$ and self-adjoint $\Gamma_0$. Then the scalar spectral measures of $\bA\ci{\Gamma_0+t\Gamma}$ are mutually singular with respect to $\nu$ for all except maybe countably many $t\in \R$.
\end{thm}

\subsection{Equivalence Classes and Spectral Multiplicity}

In \cite{GT}, Gesztesy--Tsekanovskii obtained structural results for Nevanlinna--Herglotz functions that are applicable to finite-rank perturbations. Under the assumption that ker$(I+\Gamma {\bf F}(z))=\{0\}$ for all $z\in\CC_+$,
some of these results resemble the Kato--Rosenblum Theorem \ref{t-KR} and Aronszajn--Donoghue Theorem \ref{t-AD}. 

We begin by introducing the following sets, where $1\leq r\leq d$:
\begin{align*}
    S\ti{r}\left(\bmu \right)\ti{ac}&=\left\{x\in\RR ~\Big|~ \lim_{y\to 0^+}{\bf F}(x+iy) \text{ exists finitely, and }\right.\\
    &\qquad\,\,\quad\qquad\left.\lim_{y\to 0^+}\text{rank}(\IM({\bf F}(x+i0)))=r\right\},\\
    S\left(\bmu \right)\ti{ac}&=\bigcup_{r=1}^d S\ti{r}\left(\bmu \right)\ti{ac}.
\end{align*}
Here, the existence of matrix limits are understood entrywise. Consider the equivalence classes of $S\ti{r}(\bmu\ciG)\ti{ac}$ and $S(\bmu\ciG)\ti{ac}$ associated with ${\bf F}\ciG (z)$; and denote them by $\cE_r(\bmu\ciG)\ti{ac}$ and $\cE(\bmu\ciG)\ti{ac}$, respectively.

In this setting, Gesztesy--Tsekanovskii \cite[Theorem 6.6]{GT}\footnote{Gesztesy--Tsekanovskii present these results for a slightly more general setting, when ${\bf F}$ and ${\bf F}\ciG$  are related by a certain linear fractional transformation. Their presentation reduces to ours upon making the choices $\Gamma_{1,1} = \Gamma_{2,2} = I$ and $\Gamma_{2,1} = \OZ$ and $\Gamma_{1,2} = \Gamma.$}
 have shown that:
\begin{enumerate}
    \item For $1\leq r\leq d$, the classes  $\cE_r(\bmu\ciG)\ti{ac}$, and $\cE(\bmu\ciG)\ti{ac}$ are independent of $\Gamma$.
    \item Suppose $\bmu\ci{\Gamma_1}$ is a discrete point measure for some $\Gamma_1$. Then $\bmu\ciG$ is a discrete point measure for all $\Gamma$.
    \item The set of those $x\in\RR$ for which, simultaneously, there is no $\Gamma$ such that $\lim_{y\to 0^+}\IM\left({\bf F}\ciG (x+iy)\right)$ exists and $\lim_{y\to 0^+}\det\left(\IM\left({\bf F}\ciG (x+iy)\right)\right)=0$, is a subset of $\cE_d(\bmu\ciG)\ti{ac}$.
\end{enumerate}

\section{Model Theory of Finite-Rank Unitary Perturbations}\label{s-UMT}

Taking a different route than Clark theory, we follow \cite{LT17-Clark} to set up the problem. This perspective is more natural here, since we are interested in perturbation theory. It allows us to bypass some minor technical road blocks that arise for finite-rank perturbations (when connecting the family measures with the family of operators). Some of the model theory of rank-one perturbations carries over to model theory of finite-rank setting with the added complication that one has to keep track of the order of matrix products. For example, the description of the absolutely continuous part in terms of the characteristic function has an analog for finite-rank perturbations. Other results such as identifying when the extreme situation occurs (when the de Branges--Rovnyak transcription simplifies) need to be slightly adjusted. For this particular question, taking the trace will be appropriate.

In Subsection \ref{ss-KSRKHS} we briefly mention some other representations using Krein spaces and reproducing kernel Hilbert spaces.

\subsection{Setup and Model Spaces}\label{ss-RkdModel}

Recall the setting for unitary finite-rank perturbations $\bU_\al = \bU + \bJ (\al-I)\bJ^*\bU$ with unitary $\al$, as detailed in and around \eqref{d-ual}. It is well-known that $\Ran{\bJ}$ also forms a $*$-cyclic subspace for the perturbed operators $\bU_\alpha$. Let $\bmu_\al$ be the family of matrix-valued spectral measures on $\T$ given by the Spectral Theorem through
\begin{align}\label{d-SpMeasUnit}
\bJ^*(I -  z \bU_\al^*)^{-1}\bJ
=
\int\ci{\T}\frac{d\bmu_\al(\zeta)}{1-z\bar\zeta}
\qquad\text{for }z\in \C\setminus\T.
\end{align}

It is not hard to see that the operator $\bU_\al$ is a completely non-unitary contraction for matrices $\al$ with $\|\al\|<1$. This provides us access to the associated model theory. Referring the reader to \cite[Sections 3 and 4]{LT17-Clark}, we omit the details of showing that operator $\bU_\OZ$ corresponds to the matrix-valued characteristic function
\begin{align}\label{e-teK}
\bte(z) = (K\bmu(z) - \OID)( K\bmu(z))^{-1}.
\end{align}
Here, the identity $I$ maps $\fD\to\fD$ and $K$ is the Cauchy transform of a matrix-valued measure
    $(K\bnu)(z)=\int_{\TT}\frac{d\bnu(\zeta)}{1-z \bar\zeta}.$

It is not hard to see that the relation in \eqref{e-teK} is equivalent to the Herglotz formula
\begin{align}\label{e-mHerglotz}
    (H\bmu) (z) = (\OID + \bte(z))(\OID - \bte(z))^{-1},
\end{align}
with the Herglotz transformation of a matrix-valued measure $    (H\bnu)(z)=\int_{\TT}\frac{\zeta+z}{\zeta-z}d\bnu(\zeta)
$. Now, one can reason that replacing $\bmu$ by $\bmu_\al$ in \eqref{e-mHerglotz} will result in replacing $\bte$ by $\bte\alpha^*$. And we arrive at the starting point of Aleksandrov--Clark Theory, see e.g.~\cite[Eq.~(2.5)]{Martin-uni} when $\bte(0)=\OZ$. 
 
It is worth mentioning that in starting with \eqref{d-ual} we do not really make a hidden assumption. We would recover the general starting point of \cite{Martin-uni} by taking $\bU_\al$ with strict contraction $\al$ instead of $\bU_\OZ$ with $\balpha = \OZ$.
As when dealing with rank-one perturbations, operator $\bU_\OZ$ is unitarily equivalent to the compressed shift operator on a transcription free model space.

Similar to the rank-one setting, here, the Sz.-Nagy--Foia\c s model space reduces to $H^2(\C^d)\ominus \bte H^2(\C^d)$, if and only if $\bte$ is inner (i.e.~has non-tangential boundary values that are unitary with respect to Lebesgue measure a.e.~on $\T$), if and only if $\bU$ has purely singular spectrum. See e.g.~\cite[Corollary 5.8]{LT17-Clark} for a reference of the second equivalence. Also see \cite{DL2013}.

The de Branges--Rovnyak model space reduces to one-story if and only if $\bte$ is an extreme point, and if and only if $\int\ci\T\tr (\ln(\OID - |\bte(z)|) dm(z) = -\infty$ (see \cite[Theorem 4.3.1]{Martin-uni}). There seems to be no immediate description of the extreme property in terms of the operator $\bU$ or the perturbation family $\bU\ci\al $.

In any case, the de Branges--Rovnyak model space reduces at times when Sz.-Nagy--Foia\c s model does not. When dealing with the general case of finite-rank unitary perturbations, no such reduction can be assumed a priori. This general case is the subject of Liaw--Treil \cite{LT17-Clark} and some of Martin \cite{Martin-uni} holds in this generality.

In \cite{LT17-Clark} Liaw--Treil study the general Nikolski--Vasyunin model of finite-rank Aleksandrov--Clark perturbations. Determining the unitary operator realizing this representation yields a generalization of the Clark-type operator and its adjoint. For the adjoint, the transcription choice leading to the full Sz.-Nagy--Foia\c s model features a generalization of the normalized Cauchy transform.

\subsection{Krein Spaces and Reproducing Kernel Hilbert Spaces in Applications}\label{ss-KSRKHS}

Krein spaces are indefinite inner product spaces; spaces which possess a Hermitian sesquilinear form that allows elements to have positive or negative values for their ``norm.'' A Hilbert inner product can be canonically defined on Krein spaces, so they can be viewed as a direct sum of Hilbert spaces \cite{Bog}. In particular, Krein spaces are naturally defined as extension spaces for symmetric operators with equal deficiency indices and have their own tools to determine spectral properties. Applications to the spectral analysis of direct sums of indefinite Sturm--Liouville operators is possible because so-called definitizable operators in Krein spaces are stable under finite-rank perturbations \cite{B1}. Furthermore, compact perturbations of self-adjoint operators in Krein spaces also preserve certain spectral points \cite{AJT}, and the spectral subspaces corresponding to sufficiently small surrounding neighborhoods of these points are actually Pontryagin spaces (simpler versions of Krein spaces).

Representations of symmetric operators with equal deficiency indices are also possible in reproducing kernel Hilbert spaces; Hilbert spaces of functions where point evaluation is a continuous linear functional. Among other results, Aleman--Martin--Ross \cite{AMR} carried out representations for Sturm--Liouville and Schr\"odinger (differential) operators, Toeplitz operators and infinite Jacobi matrices. The idea becomes that for each such example, the structure of the model space hosts the full information (including spectral properties) of the symmetric operator. In \cite[Section 5]{AMR}, the characteristic functions corrresponding to these examples are computed explicitly; so that the de Branges--Rovnyak model space (which is a reproducing kernel Hilbert spaces) is completely determined.

Representations in the Herglotz space constitute another interesting topic in \cite{AMR}.

\section{Spectral Theory of Finite-Rank Unitary Perturbations}\label{s-STUnit}

Consider the setting of Subsection \ref{ss-RkdModel}. Recall that $\bU_\al$ for unitary $\al$ is a unitary rank $d$ perturbation of a unitary operator $\bU$, and recall that \eqref{d-SpMeasUnit} defines the family of associated matrix-valued spectral measures $\bmu_\al$. In analogy to the self-adoint setting, we define 
the family of matrix-valued functions $W_\al$ by $d\bmu_\al(t) = W_\al (t)d\mu_\al(t)$. Taking a component-wise Randon--Nikodym derivative, we arrive at $(W_\al)\ti{ac}:= d\bmu_\al / dx$.
Further, recall that $\bte$ is the matrix-valued characteristic function of the completely non-unitary contraction $\bU_\OZ$, and that $\Delta(z) = (I-\bte^*(z)\bte(z))^{1/2}.$

A complete explicit description of the matrix-valued spectral measures of $\bU_\al $ in terms of the characteristic function is currently not available. In fact, the theory for finite-rank perturbations is lagging behind what is known for rank-one perturbations, see Theorem \ref{t-STR1Unit}. This problem has been in recent years and continues to be a field of active study. Here we explain some results in this direction.

\subsection{Spectral Properties in Terms of the Characteristic Function}

The location of the spectrum of the perturbed operator is captured by:
\begin{thm}[{see Mitkovski \cite[Corollary 4.4]{Mitkovski_IUMJ}}]
The spectrum of $\bU_\al $ consists of those points $\lambda\in \T$ at which either $\bte$ cannot be analytically continued across $\lambda$, or $\bte(\lambda)$ is analytically continuable with $\bte(\lambda)-\al $ not invertible.
\end{thm}

In combination with von Neumann's theorem, Theorem \ref{t-vN}, a characterization by Lifshitz \cite[Theorem 4]{Lifshitz} of the essential spectrum of $U_\OZ$ says that it consists of those points $\lambda\in \T$ for which (at least) one of the following conditions fails:
\begin{itemize}
\item $\bte$ is analytic on some open neighborhood of $\lambda$,
\item there is a neighborhood $N_\lambda$ of $\lambda$ so that $\bte$ is unitary for all $\lambda \in N_\lambda\cap \T$.
\end{itemize}

For the absolutely continuous part of the perturbed operator's spectral measure, a full matrix-version becomes available upon combination of Liaw--Treil \cite[Theorem 5.6]{LT17-Clark} with the Herglotz formula \eqref{e-mHerglotz} for $\bU_\al $, which is obtained from that for $\bU$ by simultaneously replacing $\bmu$ by $\bmu_\al $ and $\bte$ by $\bte\al ^*$. Namely, we have
\begin{align}\label{e-ACinTheta}
(I_n-\al  b(\lambda)^*)W_\al (\lambda) (I_n-b(\lambda)\al ^*)
=
(\Delta(\lambda))^2
\qquad\text{for Lebesgue a.e.~}\lambda \in \T,
\end{align}
in the sense of non-tangential boundary limits. 
In particular, for the absolutely continuous part, the multiplicity function is given by a non-tangential limit
rank$
\,(W_\al (\lambda))\ti{ac} = \lim_{z \to\lambda}\rk\,\Delta(z).$ Slightly weaker results are contained in Douglas--Liaw \cite{DL2013}.

\subsection{Singular Part in Terms of the Characteristic Function}

As for rank-one perturbations, capturing the singular part is a more difficult venture. The main problem here is that Poltoratski's Theorem requires a major adjustment (see Subsection \ref{ss-PoltThm} for a discussion). As a result, a description of the singular part in terms of the characteristic function is still outstanding. 

For regular points (i.e.~those that lie in the complement of the essential spectrum), both eigenvalues and eigenvectors of $\bU\ci\al $ are described in Martin \cite[Proposition 5.2.2]{Martin-uni}. 
Namely, a regular point $\lambda \in\T$ is an eigenvalue of $\bU \ci\al $ if and only if $\lim_{z\to\lambda}(\al  \bte^*(z) - \bU ^*)$ exists and is not invertible. Eigenvectors are those functions $\chi_{\{\lambda\}}{\bf x}$ with ${\bf x}\in \C^d\cap \ker (\al  \bte^*(\lambda) - \bU ^*)$ and where $\chi$ denotes the characteristic function. In that same proposition, a necessary and sufficient condition is provided for a point to not be an eigenvalue of $\bU \ci\al $ for any unitary $\al $.

There are many open questions remaining in this area. Some of them are currently being investigated.

\section{Appendix: Brief Summaries of Other Closely Related Topics}\label{s-APP}

We discuss Aleksandrov Spectral Averaging and Poltoratski's Theorem. These are both central tools in the field. Thereafter, we briefly illuminate the Simon--Wolff Theorem to which we attribute some of the popularity of the topic among mathematical physicists. We wrap up with a promising direction connecting the field to modern function theoretic operator theory.

\subsection{Aleksandrov Spectral Averaging}

Undoubtedly one of the most celebrated results of the field is the following averaging formula. On the side we mention that we can retrieve restrictions on the Aleksandrov--Clark family of spectral measures, e.g., by choosing the function $g$ to be the characteristic function of a set of Lebesgue measure zero. We being by considering the rank-one setting and will then turn to finite-rank.
For part of this subsection we follow \cite{cimaross}.

\begin{thm}[Aleksandrov dinsintegration theorem, see \cite{Alek} and {\cite[Theorem 9.4.11]{cimaross}}]\label{t-ASA}
For $g\in L^1(\T)$ we have
\begin{align}
    \int\left(\int g(\zeta)d\mu_{\al}(\zeta)\right)dm(\al)=\int g(\zeta)dm(\zeta).
\end{align}
\end{thm}

For a bounded Borel function $f$ on $\TT$, let 
\begin{align}\label{e-disintegrationtransform}
    (Gf)(\al):=\int f(\zeta)d\mu_{\al}(\zeta).
\end{align}
It is one of the main aspects of Theorem \ref{t-ASA} that for $f\in L^1(\T),$ the function $Gf$ makes sense for Lebesgue a.e.~$\alpha\in \T$ and that it is integrable.

It turns out that $G$ satisfies even more subtle mapping properties. We briefly summarize those before we explain what is known for finite-rank perturbations.

Due to the assumption that $\te(0)=0$, we see from Subsection \ref{ss-SzNF} that $\|\mu_{\al}\|=1$ and so 
\begin{align*}
\|Gf\|_{\infty}\leq \|f\|_{\infty}.    
\end{align*}
Note also that the function $Gf$ is continuous whenever $f$ is continuous. The Monotone Class Theorem (see i.e.~\cite[Theorem 9.4.3]{cimaross}) can be used to show that if $f$ is a bounded Borel function, then $Gf$ is also a bounded Borel function. Hence, the integral 
\begin{align*}
    \int_{\TT}(Gf)(\al)dm(\al),
\end{align*}
makes sense. In fact, the transformation $G$ in  \eqref{e-disintegrationtransform} can be extended to many classes of functions. Not only do we have $GC\subset C$, $CL^{\infty}\subset L^{\infty}$, and $GL^1\subset L^1$, but also $GL^p\subset L^p$ ($1\leq p\leq \infty$), $G(BMO)\subset BMO$, $G(VMO)\subset VMO$, and $GB_{pq}^s\subset B_{pq}^s$, where $B_{pq}^s$ are the Besov classes, see \cite{Alek}.

\

Now, let us turn to what is known about Aleksandrov Spectral Averaging for finite-rank perturbations.

In the unitary setting, a generalization of the Aleksandrov Spectral Averaging formula for continuous functions was obtained in Elliot \cite{Elliot2010} under extra conditions and in Martin \cite[Theorem 3.2.3]{Martin-uni}.

For self-adjoint operators a Aleksandrov-type Spectral Averaging formula was proved, Liaw--Treil \cite[Theorems 4.1, 4.6]{LTJST}. These formulas imply  restrictions on the singular parts of  families of Aleksandrov--Clark measures.

Aleksandrov Spectral Averaging for the Drury--Arveson space in the setting for inner characteristic functions was achieved by Jury \cite[Theorem 2.9]{Jury}. For more on Aleksandrov--Clark theory for the Drury--Arveson space see Subsection \ref{ss-DS} and the references therein.

\subsection{Poltoratski's Theorem}\label{ss-PoltThm} 

Deep at the heart of many results in Aleksandrov--Clark theory lies the celebrated result (proved by Poltoratski in \cite{NONTAN}) stating that for a Radon measure $\tau$ on $\T$ and $f\in L^2(\tau)$ the normalized Cauchy transform $\frac{\cC f\tau(z)}{\cC\tau(z)}$ possesses non-tangential boundary values $z\to\lambda$ for $\tau\ti{s}$-a.e.~$\lambda\in\T.$
This result is so important, because it empowers us to study the behavior of the spectral measure on sets that are of Lebesgue measure zero. In particular, one can sometimes use Poltoratski's Theorem to retrieve information about the singular parts of the spectral measures.

Direct sum examples of scalar characteristic functions immediately show that a literal extension of the statement of Poltoratski's Theorem is not possible to the finite-rank setting. Nonetheless, Kapustin--Poltoratski \cite[Theorem 3]{KP06} have proved a finite-rank analog which features a matrix-valued numerator alongside a scalar-valued denominator as well as an multiplication by a left inverse of the coordinate map $\bJ$ in \eqref{d-ual}. This left inverse `automatically' annihilates directions in which the limit of the ratio does not exist.

\subsection{Simon--Wolff Criterion}\label{ss-SIMWOL}

In \cite[Theorem 3 of Section 2]{SimonWolff} Simon and Wolff provided a characterization---formulated in terms of the spectral measure $\mu$---of when rank-one perturbation problems $A_\gamma$ are pure point for Lebesgue a.e.~parameters $\gamma\in \R.$ They applied their result to showing that the one-dimensional discrete random Schr\"odinger operator exhibits so-called Anderson localization, see \cite{JaksicLecture, SIMREV}. The idea of the Simon--Wolff localization proof was to sweep through the parameter domain for the perturbed operators' random coupling constants.

In Poltoratski \cite{Poltoratski1997} the Simon--Wolff Theorem was extended to from the rank-one to the finite-rank setting.

Simon--Wolff's celebrated work initiated further applications of rank-one perturbations to a generalization of random Schr\"odinger operators called Anderson-type Hamiltonians, see e.g.~\cite{JaksicLecture}. 
Anderson-type Hamiltonians are obtained from perturbing a self-adjoint operator by countably infintely many rank-one perturbations, each coupled by a random variable. More concretely, they are of the form
$
A_\omega=
A+
\sum \omega_i \langle\fdot, \f_i\rangle \f_i,$ where $\{\f_i\}$ forms an orthonormal basis of $\cH$, and $\omega_i$ are independent random variables that are chosen in accordance with an identical probability distribution. In view of Section \ref{s-PRELIM}, the fact that the discrete random Schr\"odinger operator features an almost surely non-compact perturbation operator underlines the level of difficulty in dealing with such objects.

In Jaksi\c c--Last \cite{JaksicLecture, JakLast2006}, these methods are utilized to prove the almost sure cyclicity of the singular spectrum of the Anderson-type Hamiltonian. And in Abakumov--Liaw--Poltoratski \cite{ALP}, it is shown that under some condition any non-trivial vector is cyclic. In Liaw \cite{L13} these results are applied to numerically support a delocalization conjecture for the 2-dimensional discrete random Schr\"odinger operator.

\subsection{Functions of Several Variables}\label{ss-DS}

Recall that the Drury--Arveson space $H^2(\mathbb{B}^n)$ is the reproducing kernel Hilbert space of functions on the open unit ball $\mathbb{B}^n$ of $\C^n$, $n\in \N$, that arises from the reproducing kernel  $k(z,w) = (1-\langle z, w\rangle\ci{\C^n})^{-1}$ with $z, w\in \mathbb{B}^n.$ 

Jury \cite{Jury} extended much of the de Branges--Rovnyak construction of Clark theory to $H^2(\mathbb{B}^n)$. As before, the de Branges--Rovnyak model spaces $\cH(\bte)$ are contractively contained in  $H^2(\mathbb{B}^n)$. The family of  Clark measure is replaced by a family of states on some noncommutative operator system. The backward shift is replaced by a canonical solution to the Gleason problem in $\cH(\bte)$. An extension of some of Jury's work to non-inner but so-called quasi-extreme characteristic functions was carried out in Jury--Martin \cite{JuryMartin1}. There, the Aleksandrov--Clark measures are necessarily generalized to certain positive linear functionals. For related work on function analytic noncommutative operator theory, we refer the reader to a series of papers by Jury and Martin \cite{JuryMartin2, JuryMartin3, JuryMartin4}.

On the side we mention that it is not immediately clear whether a perturbation problem corresponds to this Aleksandrov--Clark theory for functions of several variables.

The state of affairs for self-adjoint finite-rank perturbation problems is similar. The conditions and explicit formulas necessary to pose a well-defined problem in this area have not been investigated, to the best knowledge of the authors. However, there exists a generalization of Nevanlinna--Herglotz functions to several variables (see e.g.~\cite{LN}) whose integral representation should form a framework for the analog of the Borel transform in \eqref{e-borel}.


\subsection*{Acknowledgment}
The authors would like to thank Hari Bercovici and Tracy Weyand for pointing out references and suggesting phrases for some of the wording referring to their field of research. Special thanks to Alan Sola for reading large parts of this survey and for making useful suggestions.


\end{document}